\UseRawInputEncoding
\documentclass[12pt]{amsart}

\usepackage{hyperref}  

\usepackage{graphicx}
\usepackage{amsfonts}
\usepackage{amssymb}
\usepackage{amsmath}
\usepackage{amsthm}
\usepackage{mathrsfs}

\newtheorem{theorem}{Theorem}[section]

\newtheorem{lemma}[theorem]{Lemma}
\newtheorem{claim}{}[theorem]

\newtheorem{conjecture}[theorem]{Conjecture}

\newcommand{\del}{\backslash}

\newcommand{\cC}{\mathcal{C}}

\begin{document}

\sloppy

\title[Circuit Covers of Signed Eulerian Graphs]{Circuit Covers of Signed Eulerian Graphs}

\author{Bo Bao}
\address{Center for Discrete Mathematics, Fuzhou University, Fuzhou, P. R. China. }
\email{tomcat0830@163.com (B.Bao)}
\thanks{This research was partially supported by grants from the National Natural Sciences
Foundation of China (No. 11971111 and No. 11971110) and NSFFP (2019J01645). }

\author{Rong Chen}
\email{rongchen@fzu.edu.cn (R.Chen)}

\author{Genghua Fan}
\email{fan@fzu.edu.cn (G.Fan)}

\date{\today}

\begin{abstract}
A signed circuit cover of a signed graph is a natural analog of a circuit cover of a graph, and is equivalent to a covering of its corresponding signed-graphic matroid with circuits. It was conjectured that a signed graph whose signed-graphic matroid has no coloops has a 6-cover. In this paper, we prove that the conjecture holds for signed Eulerian graphs.
\end{abstract}

\maketitle

\section{Introduction}
Let $G$ be a graph. A {\sl signed graph} is a pair $(G,\Sigma)$ with $\Sigma\subseteq E(G)$, each edge in $\Sigma$ is labelled by $-1$ and other edges are labelled by 1. The graph $G$ can be viewed as the signed graph $(G,\emptyset)$.
A {\sl circuit} is a connected 2-regular graph.
A circuit $C$ of $G$ is \emph{balanced} if  $|C\cap\Sigma|$ is even, otherwise it is \emph{unbalanced}. We say that a subgraph of $(G,\Sigma)$ is \emph{unbalanced} if it contains an unbalanced circuit, otherwise it is \emph{balanced}. 
Signed graphs is a special class of ``biased graphs", which was defined by Zaslavsky in \cite{Zas89, Zas91}. Just as biased graphs, there are two interesting classes of matroids, the class of signed-graphic matroids and the class of even-cycle matroids, associated with signed graphs, which in fact are special classes of ``frame matroids" and ``lifted-graphic matroids" associated with biased graphs, respectively.

A {\sl barbell} is a union of two unbalanced circuits sharing exactly one vertex or a union of two vertex-disjoint unbalanced circuits together with a minimal path joining them.
A \emph{signed circuit} of $(G,\Sigma)$ is a balanced circuit or a barbell.
We say the matroid with $E(G)$ as its ground set and with the set of all signed circuits as its circuit set is the {\sl signed-graphic} matroid defined on $(G,\Sigma)$. We say that $(G,\Sigma)$ is \emph{flow-admissible} if each element of $E(G)$ is in a circuit of its signed-graphic matroid, that is, each edge of $G$ is in a signed circuit of $(G,\Sigma)$.

For a positive integer $k$, we say that a signed graph $(G,\Sigma)$ has a {\it $k$-cover} if there is a family $\cC$ of signed circuits of $(G,\Sigma)$ such that each edge of $G$ belongs to exactly $k$ members of $\cC$. For ordinary graphs $G$ (signed graph $(G,\Sigma)$ with $\Sigma =\emptyset$),  a $k$-cover of $G$ is just a family of circuits which together covers each edge of $G$ exactly $k$ times. In \cite{BJJ}, Bermond, Jackson and Jaeger proved that every bridgeless graph $G$ has a 4-cover. Fan \cite{F} proved that every bridgeless graph $G$ has a 6-cover. Together it follows that every bridgeless graph $G$ has a k-cover, for every even integer $k$ greater than 2. The only left case that $k=2$ is the famous Circuit Double Cover Conjecture: every bridgeless graph $G$ has a 2-cover, which is still open and believed to be very hard. It is somehow a surprise that it is even unknown whether there is an integer $k$ such that every signed graph $(G,\Sigma)$ has a $k$-cover.

Let $A$ and $B$ be two vertex-disjoint unbalanced circuits of length
$2m+1$. Let $G$ be the signed graph obtained from $A$ and $B$ by
joining $A$ and $B$ with two internally disjoint paths of length
$2m+1$ such that the two paths form an unbalanced circuit. Then each
signed circuit of $G$ is a barbell of $6m+3$ edges. Any $k$-cover of
$G$ contains $k|E(G)|=k(8m+4)=4k(2m+1)$ edges, which must be
divisible by $6m+3=3(2m+1)$. That is, $4k$ must be divisible by 3,
which means that $k$ cannot be 2 or 4. Thus $G$  has neither
2-covers nor 4-covers. Consider the singed graph $H$ consisting of
three unbalanced circuits of length $2m+1$ with exactly one vertex
in common. Then each signed circuit of $H$ is a barbell of $4m+2$
edges. Any $k$-cover of $H$ contains $k|E(G)|=k(6m+3)=3k(2m+1)$
edges, which must be divisible by $4m+2=2(2m+1)$. That is, $3k$ must
be divisible by 2, which means that $k$ cannot be odd. Thus $H$ has
no $k$-cover for any odd $k$. These counterexamples were first given
by Fan \cite{F-2}, who also proposed the following conjecture.

\begin{conjecture}\label{Fan}
Every flow-admissible signed graph has a $6$-cover.
\end{conjecture}

In this paper, we prove

\begin{theorem}\label{main thm}
Conjecture \ref{Fan} holds for signed Eulerian graphs.
\end{theorem}

In \cite{CLLZ}, Cheng, Lu, Luo, and Zhang proved that  each signed Eulerian graph with an even number of negative edges has a $2$-cover. We will prove Theorem \ref{main thm} from a different aspect, and our proof does not rely on their result.

This paper is organised as follows. Definitions and results needed in the proof of Theorem \ref{main thm} are given in Section 2. Theorem \ref{main thm} will be proved in Section 4 by contradiction. All ``small" signed Eulerian graphs occurring in Section 4 in the proof by contradiction are dealt with in Section 3.

\section{Preliminaries}
Let $G$ be a finite graph. Let $loops(G)$ denote the set of loops in $G$. Let $\Delta(G)$ and $\delta(G)$  be the maximal and minimal degree of $G$, respectively. For a positive integer $k$, let $V_k(G)$ be the subset of $V(G)$ consisting of degree-$k$ vertices of $G$.
A subgraph $H$ of $G$ is \emph{spanning} if $V(H)=V(G)$. In this paper, we will also use $H$ to denote its edge-set. For example, we will let $G\del H$ denote $G\del E(H)$. If exactly one component of $G$ has edges, then we say that $G$ is {\sl connected up to isolated vertices}. Evidently, a connected graph is also connected up to isolated vertices, but the converse maybe not true.

We say that $G$ is  {\sl even} if every vertex of $G$ is of even degree. If an even graph is connected, we say that it is {\sl Eulerian}. A circuit $C$ of $G$ is {\sl non-separating} if $G\del C$ is connected, otherwise, it is {\sl separating}. A \emph{theta graph} is a graph that consists of a pair of vertices joined by three internally vertex-disjoint paths.
Let $\cC$ be a circuit-decomposition of an Eulerian graph $G$. Let $H$ be a graph with $\cC$ as its vertex set, where two vertices in $H$ are adjacent if and only if their corresponding circuits in $G$ have a common vertex. We say that $H$ is {\sl determined}  by $\cC$.

\begin{lemma}\label{remove C}
Let $G$ be an Eulerian graph with $\Delta(G)\geq4$. Let $C$ be a circuit of $G$. Then there is a circuit $C^{'}$ of $G$ with $C\cap C'=\emptyset$ such that $G\backslash C^{'}$ is connected up to isolated vertices.
\end{lemma}
\begin{proof}
Since $G$ is Eulerian, $G$ has a circuit-decomposition $\cC$ containing $C$. Let $H$ be the graph determined  by $\cC$. Since $G$ is connected with $\Delta(G)\geq4$, the graph $H$ is connected with at least two vertices. Let $T$ be a spanning tree of $H$. Since $T$ has at least two degree-1 vertex, $T$ has a degree-1 vertex, say $C'$, which is not $C$. Then $C'$ is the circuit as required by the lemma.
\end{proof}

\begin{lemma}\label{2-conn}
Let $G$ be a $2$-connected graph with $\vert V(G)\vert\geq 3$. For any vertex $v$ of $G$, there is an edge $e$ of $G-v$ such that $G-V(e)$ is connected.
\end{lemma}
\begin{proof}
Let $C$ be a circuit of $G$ passing through $v$ with $|C|$ as large as possible. Evidently, $|C|\geq3$ as $\vert V(G)\vert\geq 3$ and $G$ is 2-connected. Let $e$ be an edge of $C$ that is not incident with $v$. Then $G-V(e)$ is connected, otherwise we can find a longer circuit going through $v$.
\end{proof}

A set $\Sigma'\subseteq E(G)$ is a {\sl signature} of $(G,\Sigma)$ if $(G,\Sigma)$ and $(G,\Sigma')$ have the same balanced circuits and the same unbalanced circuits. Evidently, for any edge-cut $C^*$ of $G$, the symmetric difference $\Sigma\triangle C^*$ is a signature of $(G,\Sigma)$. We say that $(G,\Sigma')$ is obtained from $(G,\Sigma)$ by \emph{switching}. The following three lemmas are well-known results on signed graphs, which will be frequently used in Section 3 without reference. Please refer to  (\cite{CDFP}, Lemma 3.5.), if the reader needs more detail about Lemma \ref{switch}.


\begin{lemma}\label{switch}
All edges of a balanced signed subgraph of $(G,\Sigma)$ can be labelled by $1$ by switching.
\end{lemma}


\begin{lemma}
Each signed theta-graph has a balanced circuit and can not have exactly two balanced circuits.
\end{lemma}

\begin{lemma}
Every $2$-edge-connected signed graph containing two edge-disjoint unbalanced circuits is flow-admissible.
\end{lemma}

In (\cite{Maca-1}, Theorem 4.2.), M\'a\v cajov\'a and \v Skoviera proved that a flow-admissible signed Eulerian graph with an odd number of negative edges contains three edge-disjoint unbalanced circuits. On the other hand, since each unbalanced Eulerian signed graph with an even number of negative edges contains two edge-disjoint unbalanced circuits, we have


\begin{lemma}\label{decomposition}
A flow-admissible unbalanced signed Eulerian graph contains two edge-disjoint unbalanced circuits.
\end{lemma}

For simplicity, we will also use $G$ to denote a signed graph defined on $G$.

\section{Signed Eulerian graphs with special circuit decompositions}
Let $k$ be a positive integer. Let $kG$ be the graph obtained from $G$ by replacing each edge in $G$ with exactly $k$ parallel edges.
Consider a graph constructed as follows. For $k\geq 3$, let $G$ be a circuit of length $k$ and $N$ be a subdivision of $2G$. Let $C$ be a circuit of $N$, we say that $C$ is {\sl small} if $|V(C)\cap V_4(N)|=2$, otherwise, $C$ is {\sl long}. When $C$ is small, we also say that each vertex in $V(C)\cap V_4(N)$ is an {\sl end} of $C$.
Let $e_1, e_2$ be edges in a small circuit of $N$ such that $\{e_1, e_2\}$ is not an edge-cut of $N$. That is, $\{e_1, e_2\}$ separates the two ends of the small circuit.
We say that the signed graph obtained from $N$ by labelling $\{e_1, e_2\}$ by $-1$ and all other edges by 1  is a {\sl necklace} of {\sl length} $k$. Evidently,  all small circuits in a necklace are balanced and all long circuits are unbalanced. Hence, the small circuits form a 1-cover in a necklace.

{\bf In the rest of this section, we will always let $G$ denote a 2-connected flow-admissible signed Eulerian graph with $\delta(G)\geq4$, and $\cC$ a circuit-decomposition of $G$, and let $H$ be the graph determined  by $\cC$}. We say that $\cC$ is {\sl optimal} if it satisfies the following properties:
\begin{itemize}
\item[(CD1)] $\cC$ is chosen with the number of unbalanced circuits as large as possible.
\item[(CD2)] subject to (CD1), $\cC$ is chosen with $|\cC|$ as large as possible.
\end{itemize}

{\bf In the rest of this section, we will always assume that $\cC$ is optimal}.
For any $C\in \cC$, we say that $C$ is a {\sl balanced} vertex of $H$ if $C$ is a balanced circuit of $G$, otherwise it is {\sl unbalanced}.

The following lemma follows immediately from (CD1), (CD2), and Lemma \ref{decomposition}.

\begin{lemma}\label{adjacent-circuits}
For every pair of adjacent vertices $C_i$ and $C_j$ in $H$, if $C_i$ is balanced, we have
\begin{enumerate}
\item $1\leq |V_G(C_i)\cap V_G(C_j)|\leq2$,
\item $C_i\cup C_j$ is balanced when $C_j$ is balanced, and
\item $C_i\cup C_j$ is not flow-admissible when $C_j$ is unbalanced.
\end{enumerate}
\end{lemma}

\begin{lemma}\label{adjacent-ub-C}
For every pair of adjacent unbalanced vertices $C_i$ and $C_j$ in $H$, if $|V_G(C_{i})\cap V_G(C_{j})|\geq 3$, then $C_{i}\cup C_{j}$ is a necklace. 
\end{lemma}
\begin{proof}
Since $C_i$ and $C_j$ are unbalanced, for any circuit decomposition $\cC'$ of $C_{i}\cup C_{j}$, either all circuits in $\cC'$ are balanced or at least two of them are unbalanced. If $C_{i}\cup C_{j}$ has an unbalanced circuit avoiding some vertex in $V_4(C_{i}\cup C_{j})$, then $C_{i}\cup C_{j}$ can be decomposed into at least three circuits and two of which are unbalanced, which is not possible as $\cC$ is optimal. So each circuit in $C_{i}\cup C_{j}$ avoiding a vertex in $V_4(C_{i}\cup C_{j})$ is balanced. Hence, $C_{i}\cup C_{j}$ is a necklace. 
\end{proof}

We say that $G$ is {\sl cover-decomposable} if $G$ can be decomposed into two proper edge-disjoint flow-admissible signed Eulerian subgraphs.

\begin{figure}[htbp]
\begin{center}
\includegraphics[height=3.5 cm]{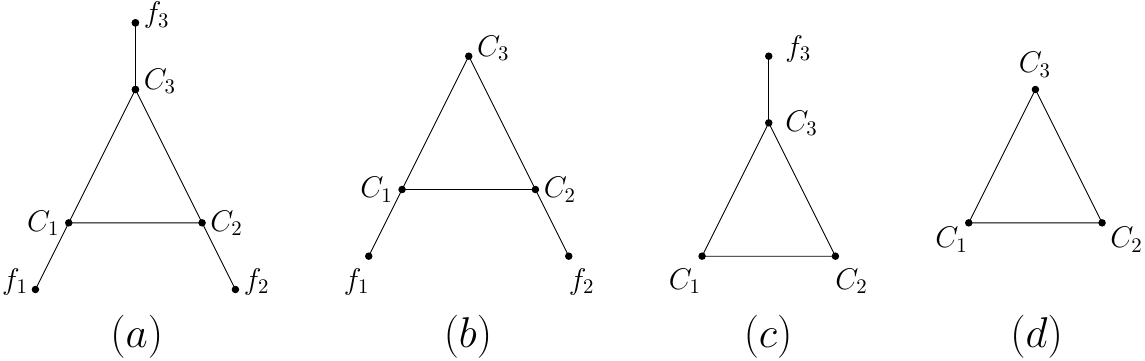}
\caption{All degree-3 vertices are balanced, and others are unbalanced. All $f_i$'s are loops of $G$. }
\label{Figure special graph}
\end{center}
\end{figure}

\begin{lemma}\label{final-structure}
If $H$ is isomorphic to a graph pictured as Figure \ref{Figure special graph} and $G$ has no balanced loops, then $G$ is cover-decomposable or has a $6$-cover.
\end{lemma}

\begin{proof}
Assume otherwise. Assume that $H$ is isomorphic to the graph pictured as Figure \ref{Figure special graph} (d). For any $1\leq i<j\leq 3$, when $|V(C_i)\cap V(C_j)|\leq 2$, it is obvious that $C_i\cup C_j$ has a $1$-cover; when $|V(C_i)\cap V(C_j)|\geq 3$,  it follows from Lemma \ref{adjacent-ub-C} that $C_i\cup C_j$ is a necklace, so $C_i\cup C_j$ has a $1$-cover too. Then $G$ has a 2-cover. So $H$ is isomorphic to a graph pictured as Figure \ref{Figure special graph} (a)-(c). Note that, $1\leq|V_G(C_i)\cap V_G(C_j)|\leq 2$ when $C_i$ is balanced by Lemma \ref{adjacent-circuits}. When some $C_i$ is a loop, implying that $G$ is isomorphic to the graph pictured as Figure \ref{Figure special graph} (b) or (c), since each pair of adjacent circuits intersect in at most 2 vertices, there are a few cases to check that $G$ has a 6-cover. So  no $C_i$ is a loop.
When $|V_G(C_i)\cap V_G(C_j)|=1$ for all $1\leq i<j\leq 3$, since $C_1\cup C_2\cup C_3$ is isomorphic to a $2K_3$-subdivision, combined the fact that all $f_i$ are unbalanced loops, it is easy to see that $G$ has a 6-cover. Hence, $|V_G(C_i)\cap V_G(C_j)|\geq2$ for some $1\leq i<j\leq 3$.

Assume that $H$ is isomorphic to a graph pictured as Figure \ref{Figure special graph} (a). By Lemma \ref{adjacent-circuits} (1) and symmetry we may assume that $V_G(C_2)\cap V_G(C_3)=\{u,v\}$. Let $C$ be the circuit of $C_2\cup C_3$ that is incident to neither $f_2$ nor $f_3$. Since $C$ is balanced by Lemma \ref{adjacent-circuits} (2), $G\del C$ is not  connected otherwise $G$ is cover-decomposable, so $V_G(C_1)\cap V_G(C_2\cup C_3)\subseteq V_G(C)-\{u,v\}$. When $|V_G(C_1)\cap V_G(C_2)|=|V_G(C_1)\cap V_G(C_3)|=1$, the graph $G$ has a 2-cover. When $|V_G(C_1)\cap V_G(C_i)|=2$ for some $2\leq i\leq3$, since $G\del loops(G)$ is balanced by Lemma \ref{adjacent-circuits} (2), there is a non-separating balanced circuit $C'$ contained in $C\cup C_i$, implying that $G$ is cover-decomposable. Hence, $H$ is isomorphic to a graph pictured as Figure \ref{Figure special graph} (b) or (c).

Assume that $V_G(C_2)\cap V_G(C_3)=\{u,v\}$. Since exactly one of $\{C_2, C_3\}$ is unbalanced, there is a $(u,v)$-path $P$ of $C_2\cup C_3$ such that a circuit in $C_2\cup C_3$ is unbalanced if and only if it contains $P$. Since all degree-3 vertices in Figure \ref{Figure special graph} are balanced, $P$ is not incident to $f_2$ or $f_3$. Let $C$ the unique balanced circuit of $C_2\cup C_3$ that is not incident to $f_2$ or $f_3$.
Since $(C_2\cup C_3)-C$ and $f_3$ are unbalanced, $G\del C$ is not connected, so $V_G(C_1)\cap V_G(C_2\cup C_3)\subseteq V_G(C)-\{u,v\}$. When $C_1\cup C_2$ is a necklace, implying that $H$ is isomorphic to a graph pictured as Figure \ref{Figure special graph} (c) by Lemma \ref{adjacent-circuits}, there is a non-separating small circuit $C'$ of the necklace $C_1\cup C_2$ with $C'\subseteq C_1\cup (C_2-P)$. Since $(C_2\cup C_3)-C$ and $f_3$ are unbalanced, $G\del C'$ is flow-admissible, so $G$ is cover-decomposable as $C'$ is balanced. Hence, by Lemma \ref{adjacent-circuits} (1) or Lemma  \ref{adjacent-ub-C}, we have $|V_G(C_1)\cap V_G(C_i)|\leq 2$ for each $2\leq i\leq 3$. Moreover, since $V_G(C_1)\cap V_G(C_2\cup C_3)\subseteq V_G(C)-\{u,v\}$, repeatedly using a similar strategy, we can find  a 6-cover of $G$ or a non-separating balanced circuit $C$ such that $G\del C$ is flow-admissible, a contradiction.

By symmetry we may therefore assume that $|V_G(C_i)\cap V_G(C_3)|=1$ for each $1\leq i\leq 2$. Set $m=|V_G(C_1)\cap V_G(C_2)|\geq2$. When $m=2$, by simple computation, the lemma holds.
So $m\geq3$. By Lemmas \ref{adjacent-circuits} and \ref{adjacent-ub-C}, $H$ is isomorphic to the graph pictured as Figure \ref{Figure special graph} (c) and $C_1\cup C_2$ is a necklace of length $m$.
Assume that $G$ is a counterexample to the lemma with $|V(G)|$ as small as possible. When $C_3$ does not share a vertex with a small circuit $C$ of $C_1\cup C_2$, delete $C$ and identify its two ends as a new vertex. Let $G'$ be the new graph. Then $G'$ is cover-decomposable or has a 6-cover by the choice of $G$, so is $G$ since $C$ is balanced. Hence, $C_3$ intersects all small circuits of $C_1\cup C_2$. Moreover, since $m\geq3$ and $|V_G(C_i)\cap V_G(C_3)|=1$ for each $1\leq i\leq 2$, there are edge-disjoint long circuits $C'_1, C'_2$ of $C_1\cup C_2$ with $|V_G(C'_1)\cap V_G(C_3)|=2$ and $|V_G(C'_2)\cap V_G(C_3)|\geq1$. Since $C'_1, C'_2$ are unbalanced and $C'_1\cup C'_2=C_1\cup C_2$, the graph determined by $\{C'_1, C'_2, C_3, \{f_3\}\}$ isomorphic to the graph pictured as Figure \ref{Figure special graph} (c). Since $|V_G(C'_1)\cap V_G(C_3)|=2$, the lemma holds by similar analysis in the second paragraph of the proof.
\end{proof}

Let $C$ be a separating circuit of a graph $G$ with $u,v\in V(C)$. Let $P$ be an $(u,v)$-path on $C$. For a component $G'$ of $G\del C$, if $V(G')\cap V(P)\neq\emptyset$ we say that $G'$ {\sl intersects} $P$; if $V(G')\cap V(C)\subseteq V(P)-\{u,v\}$ we say that $G'$ {\sl properly intersects} with $P$.

\begin{lemma}\label{cover-decomposable}
Let $C$ be a separating circuit of $G$ such that all components of $G\del C$ are unbalanced. Let $C'$ be a circuit-component of $G\del C$ with $\{u,v\}=V(C)\cap V(C')$. Let $P_1$ and $P_2$ be the $(u,v)$-paths of $C$. When $C$ is balanced or $G\del C$ has at least three components, one of the following holds.
\begin{itemize}
\item[(1)]  $G$ is cover-decomposable, or
\item[(2)] $G\del C$ has exactly three components, none of which is flow-admissible and one of which properly intersects with $P_i$ for each $1\leq i\leq2$.
\end{itemize}
\end{lemma}
\begin{proof}
 Assume that (1) is not true. Without loss of generality we may assume that $C'=\{e,f\}$.
Since $C'$ is unbalanced, we may assume that $P_1\cup\{e\}$ and $P_2\cup\{x\}$ are balanced for some $x\in\{e,f\}$. Since $G\del C$ has at least two components, besides $C'$, some component of $G\del C$  intersects with some $P_i$, say $P_2$.  Since $C$ is balanced or $G\del C$ has at least three components, $G\del (P_1\cup\{e\})$ has two edge-disjoint unbalanced circuits. Since (1) does not hold,  $G\del (P_1\cup\{e\})$ is disconnected, so there exists some non-flow-admissible component of $G\del C$ properly intersecting with $P_1$. Repeating the analysis, there is also a non-flow-admissible component of $G\del C$ properly intersecting with $P_2$. So $G\del C$ has at least three components.

Let $G_i$ be the union of the components of $G\del C$ that properly intersects with $P_i$ for each $1\leq i\leq2$. Then $G_1$ and $G_2$ are not flow-admissible. Assume that $G_1$ is disconnected. Since $G_1$ contains two edge-disjoint unbalanced circuits, $G_1\cup P_1\cup\{x\}$ and $G\del (G_1\cup P_1\cup\{x\})$ are flow-admissible, implying that (1) holds. Hence, $G_1$ is connected, so is $G_2$ by symmetry. Besides $C', G_1$ and $G_2$, assume that $G\del C$ has another component $G_3$. Since $G_3$ is unbalanced and intersects $V(P_1)$ and $V(P_2)$ by the definition of $G_1$ and $G_2$ and the fact that $G$ is 2-connected, both $G_1\cup P_1\cup \{f\}$ and $G\del (G_1\cup P_1\cup\{f\})$ are flow-admissible, a contradiction. So $G\del C$ has exactly three components $C', G_1$ and $G_2$, that is, (2) holds.
\end{proof}

For an $(u,v)$-path $P$ of $G$, we say that $P$ is  {\sl pendant}  if $u\in V_1(G)$, $v\notin V_1(G)\cup V_2(G)$ and all internal vertices of $P$ are in $V_2(G)$.

\begin{lemma}\label{star+1}
Let $H$ be a tree with a unique vertex $C$ of degree at least three, all leaf vertices are unbalanced, and all pendant paths have at most two edges. When $C$ is balanced, $V_2(H)=\emptyset$. When $C$ is unbalanced, all degree-$2$ vertices of $H$ are balanced triangles and leaf vertices that are adjacent to degree-$2$ vertices are loops. Then $G$ is cover-decomposable or has a $6$-cover.
\end{lemma}
\begin{proof}
Assume that the lemma is not true. Since $G$ has a 6-cover when each component of $G\del C$ is a loop, there is a vertex $C'$ in $H$ adjacent to $C$ with $|C'|\geq2$. Set $m=|V_G(C)\cap V_G(C')|$.
Since $G$ is 2-connected and $\delta(G)\geq4$, we have $m\geq2$.

We claim that $C'$ is balanced or $|C'|\neq2$. Assume otherwise. Then $C'$ is a component of $G\del C$ as all degree-2 vertices of $H$ are balanced. Since $|C'|=2$, we have $m=2$. Let $\{u, v\}=V_G(C')\cap V_G(C)$, $P_1$ and $P_2$ be the $(u,v)$-paths of $C$.
By Lemma \ref{cover-decomposable}, $G\del C$ has exactly three components $C', G_1$ and $G_2$, where $G_1$ and $G_2$ properly intersect $P_1$ and $P_2$, respectively. When $C\cup G_1$ is a necklace, there is a small circuit $D$ of $C\cup G_1$ such that $G\del D$ is connected. Since $C'$ and $G_2$ are unbalanced, $G\del D$ is flow-admissible, so $G$ is cover-decomposable. Hence, $G_1$ is an unbalanced circuit of size at most 2 or $G_1$ consists of a balanced triangle and a loop, so is $G_2$ by symmetry. By simple computation, $G$ is cover-decomposable or has a 6-cover.

Assume that $C'$ is balanced. Then $C'\in V_2(H)$ is a triangle. So $C$ is unbalanced and $|V_G(C')\cap V_G(C)|=2$ by Lemma \ref{adjacent-circuits}. Let $u,v, P_1, P_2$ be defined as above. Let $e$ be the loop incident with $C'$ and $f$ the edge in $C'$ whose ends are $u,v$. Since $C$ is unbalanced, $P_1\cup\{f\}$ is balanced and $P_2\cup\{f\}$ is unbalanced. Evidently, (a) a component of $G\del C$ properly intersects with $P_1$, otherwise $P_1\cup\{f\}$ and its complement are flow-admissible; and (b) no component of $G\del C$ intersects $P_2-\{u,v\}$, otherwise the union $G'$ of $P_2\cup\{f\}$ and all components of $G\del C$ intersecting $P_2-\{u,v\}$ and $G\del G'$ are flow-admissible. Then $P_2\cup (C'-\{f\})\cup\{e\}$ and its complement are flow-admissible, a contradiction.

We may therefore assume that $C'$ is unbalanced with $|C'|\geq3$, implying that $C$ is unbalanced by Lemma \ref{adjacent-circuits}. By the choice of $C'$, for each component $G'$ of $G\del C$, either $G'$ is a loop or $|G'|\geq3$. When $|G'|\geq3$, $C\cup G'$ is a necklace by Lemma \ref{adjacent-ub-C}.
Let $D$ be a small circuit of $C\cup C'$. Since $G\del D$ has two edge-disjoint unbalanced circuits, $G\del D$ is disconnected, so a component $G_D$ of $G\del C$ properly intersects in $C\cap D$. Since $C\cup C'$ has three small circuits, $G_D$ is the unique component of $G\del C$ properly intersecting in $C\cap D$ and  $C\cup C'$ has exactly three small circuits, implying $|C'|=3$, otherwise $G$ is cover-decomposable.
When $G_D$ is not a loop, there is a small circuit $D'$ of $C\cup G_D$ such that $G\del D'$ is connected, so $G$ is cover-decomposable. Hence, $G_D$ is a loop.  By the choice of $C'$, each component $G'$ of $G\del C$ that is not a loop is an unbalanced triangle. When $C'$ is the unique component of $G\del C$ that is not a loop, $G$ has a 3-cover. When there is another component $G_1$ of $G\del C$ that is not a loop, let $D$ be a small circuit of $C\cup C'$ intersecting $G_1$. Let $G'$ be the union of $D\cup G_1$ and the loop incident with $D$. Then $G'$ and $G\del G'$ are flow-admissible, so $G$ is cover-decomposable.
\end{proof}

\section{Proof of Theorem \ref{main thm}.}
In this section, we prove Theorem \ref{main thm}, which is restated here in a slightly different way.

\begin{theorem}
Every flow-admissible signed Eulerian graph has a $6$-cover.
\end{theorem}
\begin{proof}
Assume that the result is not true. Let $G$ be a counterexample with $|V(G)|$ as small as possible. Evidently,

\begin{claim}\label{4.1}
\begin{itemize}
\item $G$ is unbalanced with $\delta (G)\geq 4$;
\item $G$ has no balanced loops; and
\item $G$ is not cover-decomposable, in particular, if $C$ is a non-separating balanced circuit of $G$, then $G\del C$ is not flow-admissible.
\end{itemize}
\end{claim}

\begin{claim}\label{delete loop}
$G$ is $2$-connected.
\end{claim}
\begin{proof}[Subproof.]
Assume otherwise. There are edge-disjoint Eulerian subgraphs $G_1, G_2$ of $G$ with $|E(G_1)|, |E(G_2)|\geq2$, with $\{v\}=V(G_1)\cap V(G_2)$, and with $E(G)=E(G_1)\cup E(G_2)$. Since $G$ is a minimal counterexample and not cover-decomposable, $G_{1}$ and $G_2$ are unbalanced. Let $G_i^+$ be a signed graph obtained from $G_i$ by adding an unbalanced loop $e_i$ incident with $v$  for each integer $1\leq i\leq 2$. Since $G_1^+$ and $G_2^+$ are flow-admissible,  both of them have $6$-covers by the choice of $G$. Since $|V(G_{1})\cap V(G_{2})|=1$, we can obtain a 6-cover of $G$ by combining 6-covers of $G_1^+$ and $G_2^+$, a contradiction.
\end{proof}

Let $\cC$ be an optimal circuit decomposition of $G$ and $H$ the graph determined by $\cC$. Since $G$ is connected, so is $H$. By Lemma \ref{decomposition}, at least two members of $\cC$ are unbalanced. Hence, by Lemma \ref{adjacent-ub-C}, $|V(H)|\geq3$ and the following holds. If a block of $H$ contains exactly one cut-vertex of $H$, we say the block is a {\sl leaf block}.

\begin{claim}\label{b-circuit in H}
Each balanced vertex of $H$ is a cut-vertex, in particular, each vertex in a leaf block of $H$ that is not a cut-vertex is unbalanced.
\end{claim}

By \ref{b-circuit in H} or the third part of \ref{4.1}, for any vertex $C$ of $H$, all components of $G\del C$ are unbalanced.
For a subgraph $H'$ of $H$, 
each vertex $v \in V(H')$ is labeled by a circuit $C_v$ in $\cC$.  We say that the subgraph $G'=\cup_{v\in V(H')} E(C_v)$ corresponds to $H'$.

\begin{claim}\label{adjacent-circuits+1}
Let $e$ be a cut-edge of $H$ whose ends are $C_i$ and $C_j$. If $e$ is not a leaf edge and $H-\{C_i, C_j\}$ has exactly two components, then $C_i$ or $C_j$ is unbalanced.
\end{claim}
\begin{proof}[Subproof.]
Assume to the contrary that $C_i$ and $C_j$ are balanced. Let $G_1$ and $G_2$ be the subgraphs of $G$ corresponding to the two components of $H-\{C_i,C_j\}$ with $V(G_1)\cap V_G(C_i)\neq\emptyset$. It follows from \ref{b-circuit in H} that $G_1, G_2$ are unbalanced. Moreover, since $G$ is $2$-connected, by Lemma \ref{adjacent-circuits},  we have $|V_G(C_i)\cap V_G(C_j)|=2$.  Let $u\in V_G(G_1)\cap V_G(C_i)$ and $v\in V(G_2)\cap V_G(C_j)$. Since $ |V_G(C_i)\cap V_G(C_j)|=2$, the graph $C_i\cup C_j$ has a circuit $C$ avoiding $u$ and $v$ such that $(C_i\cup C_j)\del C$ is connected up to isolated vertices. Since $H-\{C_i, C_j\}$ has exactly two components, $G\del C$ is connected, so $G\backslash C$ is flow-admissible. Moreover, since $C_i\cup C_j$ is balanced by Lemma \ref{adjacent-circuits}, $C$ is balanced, so $G$ is cover-decomposable, a contradiction.
\end{proof}

\begin{claim}\label{remove balanced C}
For any separating circuit $C\in\cC$,  if $G'$ is a component of $G\backslash C$ that is not flow-admissible, then one of the following holds.
\begin{itemize}
\item[(1)] $G'$ is an unbalanced circuit such that $|G'|\leq2$ or $C\cup G'$ is a necklace. In particular, when $C$ is balanced, $|G'|\leq2$.
\item[(2)] $G'$ consists of a loop and a balanced triangle.
\end{itemize}
\end{claim}
\begin{proof}[Subproof.]
When $G'$ is a circuit, since $\delta(G)\geq4$, by Lemmas \ref{adjacent-circuits} and \ref{adjacent-ub-C}, (1) holds. Assume that $G'$ is not a circuit. When $G'$ consists of exactly two edge-disjoint circuits that share exactly one vertex, since $C$ only shares vertices with the balanced circuit of $G'$ and  $\delta(G)\geq4$, the unbalanced circuit $C'$ in $G'$ has at most two edges.  When $|C'|=2$, there is a non-separating balanced circuit of $G$ contained in $C\cup C'$, a contradiction. So $C'$ is a loop. By Lemma \ref{adjacent-circuits} and \ref{delete loop}, the balanced circuit in $G'$ is a triangle, so (2) holds. Hence,  we may assume that $\Delta(G')\geq6$ or $|V_4(G')|\geq2$.

Since $G'$ is not flow-admissible, by switching we may assume that there is a unique edge $e$ of $G'$ labelled by $-1$ and all other edges in $G'$ are labelled by 1.
When $e$ is a loop, let $v$ be the end of $e$, and $B$ a block of $G'\del\{e\}$ containing $v$,  and let $C'$ be a circuit of $B$ containing $v$; otherwise, let $\{v\}=\emptyset$, and $B$ the block containing $e$, and let $C'$ be a circuit of $B$ with $e\in C'$. If possible, we may further assume that $C'$ is chosen with $V_G(C')\cap V_2(G')\neq\emptyset$. By Lemma \ref{remove C}, there is a circuit $C_1$ of $G'\del loops(G')$ with $C'\cap C_1 =\emptyset$ such that $G'\del C_1$ is connected up to isolated vertices. Since $C_1$ is balanced and $G\del C_1$ has two edge-disjoint unbalanced circuits, $G\del C_1$ is not connected. Hence, $V_G(C)\cap V(G')\subseteq V_G(C_1)$ and $\emptyset\neq V_2(G')\subseteq V_G(C_1)$ as $e$ is the only edge in $G'$ which has a chance to be a loop. By the choice of $C'$,
the set $V_2(G')$ is contained in another block $B'$ of $G'$ with $B\neq B'$ as $C'$ contains no vertex in $V_2(G')$. Since $V_G(C)\cap V(G')\subseteq V(B')$ and $G$ is $2$-connected, $|B|=1$, a contradiction to the choice of $B$.
\end{proof}

\begin{claim}\label{cn after delete C}
For any $C\in\cC$, the graph $G\del C$ has at most two components.
\end{claim}
\begin{proof}[Subproof.]
Assume that $G\del C$ has three components. Since each component $G'$ of $G\del C$ is unbalanced and $G\del G'$ is flow-admissible, $G'$ is not flow-admissible . By \ref{remove balanced C}, $H$ is a tree with $C$ as a unique vertex of degree at least three, and all its pendant paths have at most two edges.
When $C$ is balanced, \ref{adjacent-circuits+1} implies that $V_2(H)=\emptyset$. Hence, by \ref{remove balanced C} and Lemma \ref{star+1}, $G$ is cover-decomposable or has a 6-cover, a contradiction.
\end{proof}

\begin{claim}\label{C&C'+}
For any balanced vertex $C$ of $H$, each degree-$1$ vertex of $H$ adjacent with $C$ is a loop of $G$.
\end{claim}
\begin{proof}[Subproof.]
Let $C'$ be a degree-$1$ vertex of $H$ adjacent with $C$. Assume that $C'$ is not a loop of $G$. Then $|C'|=|V_G(C)\cap V_G(C')|=2$ by \ref{remove balanced C}. It follows from Lemma \ref{cover-decomposable} and \ref{cn after delete C} that $G$ is cover-decomposable, a contradiction.
\end{proof}

\begin{claim}\label{tree-case}
$H$ is not a tree.
\end{claim}
\begin{proof}[Subproof.]
Assume otherwise. By \ref{cn after delete C}, $H$ is a path. Evidently, at most one vertex in $V_2(H)$ is unbalanced, otherwise, $G$ is cover-decomposable.  By \ref{adjacent-circuits+1}, no balanced vertices of $H$ are adjacent, so $|V(H)|\leq 5$. Moreover, if $|V(H)|\geq4$, then exactly one vertex in $V_2(H)$ is unbalanced.
Assume that $H$ has two adjacent vertices $C_1, C_2$ with $|V_G(C_1)\cap V_G(C_2)|\geq3$. Then $C_1\in V_1(H)$, $|V(H)|\leq4$ and $C_1\cup C_2$ is a necklace by Lemma \ref{adjacent-ub-C}. Let $C_3$ be the other vertex adjacent to $C_2$ in $H$. When $V_G(C_2)\cap V_G(C_3)$ is in a small circuit of $C_1\cup C_2$, the graph $G$ has a 6-cover.
When $V_G(C_2)\cap V_G(C_3)$ is not in a small circuit of $C_1\cup C_2$, implying $|V_G(C_2)\cap V_G(C_3)|=2$, since $V_G(C_1)\cap V_G(C_3)=\emptyset$, the graph $C_1\cup C_2$ can be decomposed to two long circuits $C'_1, C'_2$, where both share exactly one vertex with $C_3$.
Note that the circuit decomposition $(\cC-\{C_1, C_2\})\cup\{C'_1, C'_2\}$ is still optimal.
Hence, the graph determined by $(\cC-\{C_1, C_2\})\cup\{C'_1, C'_2\}$ is isomorphic to a graph pictured as Figure \ref{Figure special graph} (c) or (d). Lemma \ref{final-structure} implies that $G$ is cover-decomposable or has a 6-cover. Therefore, combined with Lemma \ref{adjacent-circuits} we can assume that every pair of adjacent vertices in $H$ share at most two vertices in $G$. Note that each degree-1 vertex of $H$ adjacent to a balanced vertex is a loop by \ref{C&C'+}. By simple computation, $G$ has a 6-cover, a contradiction.
\end{proof}

\begin{claim}\label{leaf-block}
$H$ is not $2$-connected and whose leaf blocks are isomorphic to $K_{2}$.
\end{claim}
\begin{proof}[Subproof.]
Assume otherwise. When $H$ is not 2-connected, let $B$ be a leaf block of $H$ that is not isomorphic to $K_2$, and $v$ be the unique cut-vertex of $H$ in $V(B)$. When $H$ is 2-connected, let $B=H$ and $v$ any vertex of $B$.
By Lemma \ref{2-conn}, there is an edge $e$ in $B-v$ such that $B-V_H(e)$ is connected, so $H-V_H(e)$ is also connected. Without loss of generality assume that $C_{1}$ and $C_{2}$ are the ends of $e$. Then $C_{1}\cup C_{2}$ and $G\del (C_{1}\cup C_{2})$ are connected.
Since $C_{1}\cup C_{2}$ is flow-admissible by \ref{b-circuit in H}, the graph $G\del (C_{1}\cup C_{2})$ is not flow-admissible. Since $H$ is not isomorphic to the graph pictured as Figure \ref{Figure special graph} (d) by Lemma \ref{final-structure}, $H$ has exactly three unbalanced vertices and exactly two leaf blocks, one of which is $B$ that is isomorphic to $K_3$ and the other is isomorphic to $K_2$. Let $C_{1}C_{2}C_{3}\ldots C_{n}$ be a longest path in $H$. It follows from \ref{adjacent-circuits+1} that $n=4$. By \ref{C&C'+}, the circuit $C_4$ is a loop of $G$. That is, $H$ is isomorphic to the graph pictured as Figure \ref{Figure special graph} (c). Hence, $G$ is cover-decomposable or has a 6-cover by Lemma \ref{final-structure}, a contradiction.
\end{proof}

Let $B$ be a block of $H$ with $|V(B)|\geq3$. By \ref{tree-case} and \ref{leaf-block}, such $B$ exists and $B$ is not a leaf block. When $H$ has two blocks that are not isomorphic to $K_2$, it follows from \ref{b-circuit in H} and \ref{leaf-block} that $G$ is cover-decomposable.
Hence, $B$ is the unique block of $H$ that is not isomorphic to $K_2$. By \ref{b-circuit in H}, each vertex in $B$ that is not a cut-vertex of $H$ is unbalanced.

Let $u\in V(B)$ be a cut vertex of $H$. When $u$ is unbalanced or $H$ has two pendant paths using $u$, let $H_1$ be the union of all pendant paths containing $u$, and $G_1$ the subgraph of $G$ corresponding to $H_1$. Since $|V(B)|\geq3$, by \ref{b-circuit in H} and \ref{leaf-block}, both $G_1$ and $G\del G_1$ are flow-admissible, a contradiction. Hence, $u$ is balanced and $H$ has exactly one pendant path using $u$. By the arbitrary choice of $u$, all cut-vertices of $H$ in $B$ are balanced.  Using a similar strategy, all vertices in $V_2(H)-V(B)$ are balanced. Combined with \ref{adjacent-circuits+1}, we have $V_2(H)-V(B)=\emptyset$. That is, each pendant path of $H$ has exactly one edge. By \ref{C&C'+}, each vertex in $V_1(H)$ is a loop of $G$.

When there is a vertex in $V(B)$ that is not a cut-vertex of $H$, let $v$ denote such a vertex. Otherwise, let $v$ be any vertex of $B$. By Lemma \ref{2-conn}, there is an edge $e\in B-v$ such that $B-V(e)$ is connected. Let $H_1$ be the union of $e$ and all pendant paths of $H$ using an end of $e$, and $G_1$ be the subgraph of $G$ corresponding to $H_1$. Since each vertex in $B$ that is not a cut-vertex of $H$ is unbalanced, $H_1$ contains two unbalanced vertices, so $G_1$ is flow-admissible. Since $H-V(H_1)$ is connected and has an unbalanced vertex, $H$ is isomorphic to a graph pictured as Figure \ref{Figure special graph} (a) or (b). Lemma \ref{final-structure} implies that $G$ is cover-decomposable or has a 6-cover, a contradiction.
\end{proof}

\section*{Acknowledgements} We thank the referees for their careful reading of this
manuscript and their detailed comments.


\begin{thebibliography}{}
\bibitem{BJJ}
J. C. Bermond, B. Jackson, F. Jaeger, Shortest covering of graphs with cycles, J. Combin. Theory Ser. B, 35 (1983) 297-308.

\bibitem{CDFP}
R. Chen, M. DeVos, D. Funk, I. Pivotto, Graphic representation of graphic frame matroids, Graphs and Combin., 31 (2015), 2075-2086.

\bibitem{CLLZ}
J. Cheng, Y. Lu, R. Luo, C. Q. Zhang, Shortest circuit cover of signed graphs, J. Combin. Theory Ser. B, 134 (2019), 164-178.

\bibitem{F}
G. Fan, Integer flows and cycle covers, J. Combin. Theory Ser. B, 54 (1992) 113-122.

\bibitem{F-2}
G. Fan, Flows and circuit covers in signed Graphs, Lectures in NSFC Tianyuan Summer School, Jinhua, 2018.

\bibitem{Maca-1}
E. M\'a\v cajov\'a, M. \v Skoviera, Nowhere-zero flows on signed Eulerian graphs, SIAM J. Discrete Math., 31 (2017),  1937-1952.


\bibitem{Zas89}
T. Zaslavsky, Biased graphs. I. Bias, balanced, and gains, J. Combin. Theory Ser. B, 47 (1989),  32-52.

\bibitem{Zas91}
T. Zaslavsky, Biased graphs. II. The three matroids. J. Combin. Theory Ser. B, 51  (1991),  46-72.


\end{thebibliography}
\end{document}